\documentclass[english,letterpaper,11pt,reqno]{amsart}

\usepackage{appendix}
\usepackage{amsbsy}
\usepackage{amsfonts}
\usepackage{amsmath}
\usepackage{amssymb}
\usepackage{amsthm}
\usepackage{graphicx}
\usepackage{ifthen}
\usepackage{textcomp}
\usepackage{soul}
\usepackage{enumitem,kantlipsum}
\usepackage{tikz}
\usepackage{cancel}

\usepackage{dsfont}
\usepackage[bookmarksnumbered,colorlinks]{hyperref}
\hypersetup{colorlinks=true, linkcolor=blue, citecolor=blue}
\newcommand{\lra}{\longrightarrow}

\newcommand{\RR}{\mathbb{R}}

\newtheorem*{definition-non}{Definition}
\newtheorem*{theorem-non}{Theorem}
\newtheorem*{proposition-non}{Proposition}
\newtheorem*{lemma-non}{Lemma}
\newtheorem*{corollary-non}{Corollary}

\newcommand{\beqa}{\begin{eqnarray}}
\newcommand{\beq}{\begin{equation}}
\newcommand{\eeqa}{\end{eqnarray}}
\newcommand{\eeq}{\end{equation}}

\newcommand\imp{\hspace{.2in}\Rightarrow\hspace{.2in}}

\newcommand\kk{T}

\newcommand\cd[2]{\nabla_{\!#1}{#2}}
\newcommand\cds[2]{\nabla^{\scalebox{0.4}{\emph{L}}}_{\!#1}{#2}}

\newcommand\gL{g_{\scalebox{0.4}{\emph{L}}}}
\newcommand\comma{\hspace{.2in},\hspace{.2in}}

\makeatletter
\newcommand*\owedge{\mathpalette\@owedge\relax}
\newcommand*\@owedge[1]{%
  \mathbin{%
    \ooalign{%
      $#1\m@th\bigcirc$\cr
      \hidewidth$#1\m@th\wedge$\hidewidth\cr
    }%
  }%
}
\makeatother

\makeatletter
\newcommand*{\defeq}{\mathrel{\rlap{%
                     \raisebox{0.24ex}{$\m@th\cdot$}}%
                     \raisebox{-0.24ex}{$\m@th\cdot$}}%
                     =}
\makeatother

\newcommand{\kn}{\ensuremath{\raisebox{.04 em}{\,${\scriptstyle \owedge}$\,}}}  

\usepackage{enumitem, color, amssymb}

\begin{document}
\title[]{Riemannian counterparts to Lorentzian space forms}
\author[]{Amir Babak Aazami}
\address{Clark University\hfill\break\indent
Worcester, MA 01610}
\email{Aaazami@clarku.edu}

\maketitle
\begin{abstract}
On a smooth $n$-manifold $M$ with $n \geq 3$, we study pairs $(g,T)$ consisting of a Riemannian metric $g$ and a unit length vector field $T$ with geodesic flow and integrable normal bundle.  Motivated by how Ricci solitons generalize Einstein metrics via a distinguished vector field, we propose to generalize space forms by considering those pairs $(g,T)$ whose corresponding \emph{Lorentzian} metric $\gL \defeq g - 2T^{\flat} \otimes T^{\flat}$ has constant curvature. We show by examples that  such pairs exist when $M$ is noncompact, and that complete metrics exist among them.  When $M$ is compact, however, the situation is more rigid.  In the compact setting, we prove that the only pairs $(g,T)$ whose corresponding Lorentzian metric $\gL$ is a space form are those where $(M,g)$ is flat and its universal covering splits isometrically as a product $\RR \times N$.  The nonexistence of compact Lorentzian spherical space forms plays a key role in our proof.
\end{abstract}

\section{Introduction}
\label{sec:Intro}
In this paper we search for generalizations of space forms, taking as our motivation the way in which Ricci solitons generalize Einstein metrics\,---\,namely, via a distinguished vector field.  In our case, however, the vector field arises as a ``bridge" between Riemannian and Lorentzian geometry: when a nowhere vanishing vector field $T$ is present, any Riemannian metric $g$ has a Lorentzian sibling
\beqa
\label{eqn:LS0}
\gL \defeq g - 2T^{\flat} \otimes T^{\flat},
\eeqa
and vice-versa: $g = \gL + 2T^{\flat_{\scriptscriptstyle L}} \otimes T^{\flat_{\scriptscriptstyle L}}$ (note that $T^{\flat_{\scriptscriptstyle L}} = -T^{\flat} $, and that we are assuming for convenience here that $T$ has unit length).  In this paper we therefore adopt the following strategy: we seek to classify those Riemannian manifolds whose \emph{Lorentzian} sibling \eqref{eqn:LS0} has constant curvature.  Of course, the relationship \eqref{eqn:LS0} is well known and has been studied extensively; see, e.g., \cite{olea} for a recent analysis which includes, among other things, curvature formulae.  Generally speaking, \eqref{eqn:LS0} can be a fruitful method by which to construct distinguished examples of Riemannian or Lorentzian manifolds. That is because in certain cases it allows properties of one metric to be more or less directly inferred from those of the other; e.g., distinguished curvature or geodesic completeness.  Indeed, a well known instance of the latter is the following: if $T^{\flat}$ is bounded on $TM$, then $\gL$ will be complete if $g$ is (this is a consequence of \cite[Proposition 3.4]{candela}).  Of course, one expects that the more conditions one places on $T$, the easier to read these shared properties become\,---\,the easier it becomes, so to speak, to cross the bridge ``with nice properties in hand."  Therefore we have imposed, in addition to $T$ having unit length, two further properties on $T$, properties which we motivate below and which are borne out by examples.  What we are then able to achieve in our Theorem is a classification when the underlying manifold is \emph{compact}, a classification of those Riemannian pairs $(g,T)$ whose Lorentzian counterparts are precisely the space forms:
\begin{theorem-non}
\label{thm:1}
Let $M$ be an $n$-manifold $(n \geq 3)$ and $(g,T)$ a Riemannian metric $g$ on $M$ and a unit length vector field $T$ with geodesic flow and integrable normal bundle.  Then the Riemann curvature 4-tensor of $(M,g)$ is
\beqa
\label{eqn:1}
\emph{\text{Rm}} = \frac{1}{2}\lambda g \kn g - 2\lambda g \kn (T^{\flat}\otimes T^{\flat}) - \nabla T^{\flat} \kn \nabla T^{\flat},
\eeqa
where $\kn\!$ is the Kulkarni-Nomizu product and $\lambda$ is a constant, if and only if the corresponding Lorentzian metric $g_{\scalebox{0.4}{L}}\!\defeq\!g - 2T^{\flat} \otimes T^{\flat}$ has constant curvature $\lambda$.  If $M$ is compact, then $\lambda = 0$, $(M,g)$ is flat, and its universal covering splits isometrically as a product $\RR \times N$.
\end{theorem-non}

If $M$ is simply connected (and noncompact), then $T$ is a gradient and $\nabla T^{\flat}$ is its Hessian; generally speaking, $T$ is always at least locally a gradient.  Observe how \eqref{eqn:1} generalizes constant curvature, which is the first term $\frac{1}{2}\lambda g \kn g$; Remark 1 of Section \ref{sec:conc} will further elaborate on this.  The crucial fact, however, is that the corresponding Lorentzian metric $\gL$ has constant curvature $\lambda$\,---\,and this plays the key role in our proof of the case when $\lambda > 0$.  Indeed, by \cite{CM} and \cite{Klingler}, there are \emph{no} such Lorentzian metrics when $\lambda  >0$ and $M$ is compact, a foundational result in Lorentzian geometry (see \cite{lundberg} for a comprehensive account).  As for the two properties imposed on our unit length vector field $T$\,---\,namely, that it have geodesic flow and integrable normal bundle\,---\,they are shared across $g$ and $\gL$ and ensure that the 2-tensor $\nabla T^{\flat}$ is \emph{symmetric}; indeed, the symmetry of $\nabla T^{\flat}$ is equivalent to these two properties.  By virtue of this symmetry, the curvature 4-tensors of $g$ and $\gL$, denoted by $\text{Rm}$ and $\text{Rm}_{\scalebox{0.4}{\emph{L}}}$, respectively, have a very simple relationship to each other:
\beqa
\label{eqn:rel1}
\text{Rm}_{\scalebox{0.4}{\emph{L}}} = \text{Rm} + \nabla T^{\flat} \kn \nabla T^{\flat}.
\eeqa
This will no longer hold if either of these two conditions on $T$ is dropped.  \emph{Indeed, another motivation for these properties of $T$ is that such a vector field exists prominently in \emph{de Sitter spacetime}, the canonical (and noncompact) Lorentzian manifold of constant positive curvature}; see Example 3 in Section \ref{sec:ex}.  We close this Introduction with three more remarks.  First, Riemannian manifolds satisfying \eqref{eqn:1} certainly do exist in the \emph{non}compact setting. We furnish two such examples in Section \ref{sec:Prop}, both with $\lambda > 0$; note that the first of these is complete and exists in all dimensions $\geq 3$. (This makes it clear that when $\lambda > 0$, the obstruction arises from compactness, not completeness.)  Second, in the proof of our Theorem, while the case $\lambda > 0$ makes direct use of \eqref{eqn:rel1} and relies on \cite{CM} and \cite{Klingler}, the case $\lambda \leq 0$ relies instead on a Bochner technique; see \eqref{eqn:Bochner1} in Section \ref{sec:proof}. Finally, compact Lorentzian manifolds with constant \emph{negative} curvature certainly do exist; examples can be found, e.g., in \cite{kulkarni} and \cite{Goldman}.  So do flat ones: any $(\mathbb{S}^1 \!\times\! N,-dt^2\oplus h)$, with $(N,h)$ a compact flat Riemannian manifold, yields an example, with $\nabla^{\scalebox{0.4}{\emph{L}}}t$ serving the role of $T$ above.

\section{The Riemannian-to-Lorentzian correspondence}
\label{sec:Prop}
The proof of our Theorem is a corollary of the following more general fact:

\begin{proposition-non}
Let $(M,g)$ be a Riemmanian $n$-manifold equipped with a unit length vector field $T$ with geodesic flow and integrable normal bundle. Then the Lorentzian metric $\gL$ defined by
$$
\gL \defeq g - 2 T^{\flat} \otimes T^{\flat}
$$
has Riemann curvature 4-tensor $\emph{\text{Rm}}_{\scalebox{0.4}{L}}$ given by
$$
\emph{\text{Rm}}_{\scalebox{0.4}{L}} = \emph{\text{Rm}} + \nabla T^{\flat} \kn \nabla T^{\flat},
$$
where \emph{$\text{Rm}$} is the Riemann curvature 4-tensor of $g$.
\end{proposition-non}

\begin{proof}
For $\text{Rm}$, we adopt the sign convention
$$
\text{Rm}(a,b,c,d) = g(\cd{a}{\!\cd{b}{\,c}},d) - g(\cd{b}{\cd{a}{\,c}},d) - g(\cd{[a,b]}{c},d),
$$
and similarly for $\text{Rm}_{\scalebox{0.4}{\emph{L}}}$.  With respect to the Lorentzian metric $\gL$, $T$ is now unit length ``timelike": $\gL(T,T) = -1$.  Next, denoting by $\nabla^{\scalebox{0.4}{\emph{L}}}$ the Levi-Civita connection of $\gL$, the Koszul formula shows that $T$ will have $\gL$-geodesic flow, since
$$
\cds{T}{T} = -\cd{T}{T} = 0,
$$
where $\nabla$ is the Levi-Civita connection of $g$.  Furthermore, it's clear from the  definition of $\gL$ that the $\gL$-normal bundle $T^{\perp_{\tiny L}} \subseteq TM$ will remain integrable (and equal to the $g$-normal bundle $T^{\perp}$), so that the endomorphism
\beqa
\label{eqn:D0}
\mathcal{D}_{\scalebox{0.4}{\emph{L}}} \colon T^{\perp_{\tiny L}} \lra T^{\perp_{\scriptscriptstyle L}}  \comma X \mapsto \cds{X}{T}
\eeqa
is self-adjoint with respect to the (positive-definite) metric $\gL|_{T^{\perp_{\tiny L}}} = g|_{T^{\perp}}$ induced on $T^{\perp_{\tiny L}}$.  Thus we have, in a neighborhood of each point of $M$, an orthonormal basis of eigenvectors $\{X_1,\dots,X_{n-1}\} \subseteq T^{\perp_{\tiny L}}$ of $\mathcal{D}_{\scalebox{0.4}{\emph{L}}}$,
$$
\cds{X_i}{T} = \lambda_i X_i \comma i = 1,\dots,n-1,
$$
where the eigenvalues $\lambda_i$ are smooth functions defined on said neighborhoods. In fact these are the same eigenvalues and eigenvectors as those of the endomorphism
\beqa
\label{eqn:D1}
\mathcal{D} \colon T^{\perp} \lra T^{\perp}  \comma X \mapsto \cd{X}{T},
\eeqa
since $\cds{X_i}{T} = \cd{X_i}{T}$ via the Koszul formula; in particular, $\cds{T}{X_i} = \cd{T}{X_i}$, and the remaining covariant derivatives are, by the Koszul formula again,
$$
\cds{X_i}{X_j} = 2\lambda_i\delta_{ij}T + \cd{X_i}{X_j}.
$$
In fact \eqref{eqn:D1} extends trivially to all of $TM$; furthermore, because $\cd{T}{T} = 0$ and $g(T,T) = 1$, this extension remains self-adjoint with respect to $g$, thus defining a symmetric 2-tensor which is precisely $\nabla T^{\flat}$ mentioned above (although we won't need this explicitly, observe in passing that $\nabla^{\scalebox{0.4}{\emph{L}}} T^{\flat_{\scriptscriptstyle L}} = \nabla T^{\flat}$).  Now we compute the components of the curvature tensor $\text{Rm}_{\scalebox{0.4}{\emph{L}}}$ with respect to the frame of eigenvectors above:
$$
\{T,X_1,\dots,X_{n-1}\} \comma \cds{X_i}{T} = \lambda_i X_i = \cd{X_i}{T}.
$$
Note that this frame is both $g$- and $\gL$-orthonormal, but with the important difference  that
$$
g(T,T) = 1 = -\gL(T,T).
$$
As mentioned in the Introduction, more general curvature formulae are provided by \cite{olea}, in particular of the components \eqref{eqn:Ric1} and \eqref{eqn:Ric2} appearing below; however, the most important curvature component in our Proposition, namely \eqref{eqn:Ric4} below, has not (to the best of our knowledge) appeared explicitly in the literature.  As it is also the most technical, we have therefore chosen to compute \emph{all} the curvature components explicitly.
 
\begin{enumerate}[leftmargin=*]
\item[1.] We start with the components 
$\text{Rm}_{\scalebox{0.4}{\emph{L}}}(X_i,T,T,X_j)$:
\beqa
\text{Rm}_{\scalebox{0.4}{\emph{L}}}(X_i,T,T,X_j) \!\!&=&\!\! \gL(\cds{X_i}{\cancelto{0}{\!\cds{T}{T}}},X_j) - \!\!\!\!\!\underbrace{\,\gL(\cds{T}{\overbrace{\,\cds{X_i}{T}\,}^{\lambda_i X_i}},X_j)\,}_{\text{$T(\lambda_i)\delta_{ij}+\lambda_i\gL(\cds{T}{X_i},\,X_j)$}}\!\!\!\! - \underbrace{\,\gL(\cds{[X_i,T]}{T},X_j)\,}_{\text{$\gL(\cds{X_j}{T},[X_i,T])$}}\nonumber\\
&=& -T(\lambda_i)\delta_{ij} - \lambda_i\gL(\cds{T}{X_i},X_j) - \lambda_j\!\!\!\!\!\!\!\!\!\!\!\!\underbrace{\,\gL(X_j,[X_i,T])\,}_{\text{$\gL(X_j,\cds{X_i}{T})-\gL(X_j,\cds{T}{X_i})$}}\nonumber\\
&=& -\big(T(\lambda_i) + \lambda_j\lambda_i\big)\delta_{ij} + \underbrace{\,\gL(\cds{T}{X_i},X_j)\,}_{\text{$g(\cd{T}{X_i},X_j)$}}(\lambda_j-\lambda_i)\label{eqn:Bochner}\\
&=& \text{Rm}(X_i,T,T,X_j),\label{eqn:Ric1}
\eeqa
where the last equality is due to the fact that the same answer would have been reached starting with
$$
\text{Rm}(X_i,T,T,X_j) = g(\cd{X_i}{\!\cd{T}{T}},X_j) - g(\cd{T}{\cd{X_i}{T}},X_j) - g(\cd{[X_i,T]}{T},X_j).
$$
It is worthwhile to pause here and see what would have resulted had $\cd{T}{T} \neq 0$. Other than the fact that the Lie brackets $[X_i,T]$ are $g$- and $\gL$-orthogonal to $T$ if and only if $\cd{T}{T} = 0$ (which allowed us to simplify the last term),
the first term of $\text{Rm}_{\scalebox{0.4}{\emph{L}}}(X_i,T,T,X_j)$ would have been
\beqa
\gL(\cds{X_i}{\!\cds{T}{T}},X_j) \!\!&=&\!\! X_i\big(\underbrace{\gL(\cds{T}{T},X_j)}_{\text{$-g(\cd{T}{T},X_j)$}}\big) -\!\!\!\!\!\!\! \underbrace{\,\gL(\cds{T}{T},\cds{X_i}{X_j})\,}_{\text{$-g(\cd{T}{T},2\lambda_i\delta_{ij}T+\cd{X_i}{X_j})$}}\nonumber\\
&=&\!\! -g(\cd{X_i}{\!\cd{T}{T}},X_j) +2\lambda_i\delta_{ij} \cancelto{0}{g(\cd{T}{T},T)}.\nonumber
\eeqa

\item[2.] Next, the components $\text{Rm}_{\scalebox{0.4}{\emph{L}}}(X_i,X_k,T,X_j)$, with $k \neq i$:
\beqa
\text{Rm}_{\scalebox{0.4}{\emph{L}}}(X_i,X_k,T,X_j) \!\!&=&\!\! \underbrace{\,\gL(\cds{X_i}{{\overbrace{\,\cds{X_k}{T}\,}^{\lambda_k X_k}}},X_j)\,}_{g(\cd{X_i}{\!\cd{X_k}{T}},\,X_j)} - \underbrace{\,\gL(\cds{X_k}{\overbrace{\,\cds{X_i}{T}\,}^{\lambda_i X_i}},X_j)\,}_{\text{$g(\cd{X_k}{\!\cd{X_i}{T}},X_j)$}}\nonumber\\
&&\hspace{1.75in} - \underbrace{\,\gL(\cds{[X_i,X_k]}{T},X_j)\,}_{\text{$\gL(\cds{X_j}{T},[X_i,X_k])$}}\nonumber\\
&=&\!\! g(\cd{X_i}{\!\cd{X_k}{T}},X_j) - g(\cd{X_k}{\!\cd{X_i}{T}},X_j) - \underbrace{\,g(\cd{X_j}{T},[X_i,X_k])\,}_{\text{$g(\cd{[X_i,X_k]}{T},\,X_j)$}}\nonumber\\
&=&\!\! \text{Rm}(X_i,X_k,T,X_j),\label{eqn:Ric2}
\eeqa
where in the first and second equalities we've used the fact that
$
\gL(\cdot,X) = g(\cdot,X)
$
for any $X \in T^{\perp}$, which includes $[X_i,X_k]$ by integrability of $T^{\perp}$. 
\item[3.] Finally, the components $\text{Rm}_{\scalebox{0.4}{\emph{L}}}(X_k,X_i,X_j,X_l)$, with $k \neq i$, $l \neq j$:
\beqa
\text{Rm}_{\scalebox{0.4}{\emph{L}}}(X_k,X_i,X_j,X_l) \!\!&=&\!\! \underbrace{\,\gL(\cds{X_k}{\!\cds{X_i}{X_j}},X_l)\,}_{(a)} - \underbrace{\,\gL(\cds{X_i}{\!\cds{X_k}{X_j}},X_l)\,}_{(b)}\nonumber\\
&&\hspace{1.75in} - \underbrace{\,\gL(\cds{[X_k,X_i]}{X_j},X_l)\,}_{(c)},\nonumber
\eeqa
will be computed term-by-term, using the Koszul formula, beginning with $(a)$:
\beqa
\gL(\cds{X_k}{\!\cds{X_i}{X_j}},X_l) \!\!&=&\!\! \gL(\cds{X_k}{\underbrace{\cd{X_i}{X_j}}_{``Z\,"}},X_l) + 2\lambda_i\delta_{ij}\,\gL(\underbrace{\cds{X_k}{T}}_{\lambda_k X_k},X_l)\nonumber\\
&=&\!\!\frac{1}{2}\Big[X_k(\gL(Z,X_l)) + Z(\gL(X_l,X_k)) - X_l(\gL(X_k,Z))\nonumber\\
&&\hspace{.2in}-\gL(Z,[X_k,X_l]) - \gL(X_l,[Z,X_k]) + \gL(X_k,[X_l,Z])\Big]\nonumber\\
&&\hspace{2in}+\, 2\lambda_i\lambda_k\delta_{ij}\delta_{kl}.\nonumber
\eeqa
Every ``$\gL(\cdot,\cdot)$" here contains at least one $X \in T^{\perp}$, so that it equals $g(\cdot,\cdot)$.  We thus have that
\beqa
\gL(\cds{X_k}{\!\cds{X_i}{X_j}},X_l) \!\!&=&\!\! \frac{1}{2}\Big[X_k(g(Z,X_l)) + Z(g(X_l,X_k)) - X_l(g(X_k,Z))\nonumber\\
&&\hspace{.2in}-g(Z,[X_k,X_l]) - g(X_l,[Z,X_k]) + g(X_k,[X_l,Z])\Big]\nonumber\\
&&\hspace{2in}+\, 2\lambda_i\lambda_k\delta_{ij}\delta_{kl}\nonumber\\
&=&\!\! g(\cd{X_k}{Z},X_l) + 2\lambda_i\lambda_k\delta_{ij}\delta_{kl}\nonumber\\
&=&\!\! g(\cd{X_k}{\!\cd{X_i}{X_j}},X_l) + 2\lambda_i\lambda_k\delta_{ij}\delta_{kl}.\label{eqn:term1}
\eeqa
Similarly with $(b)$,
\beqa
\gL(\cds{X_i}{\!\cds{X_k}{X_j}},X_l) \!\!&=&\!\! \gL(\cds{X_i}{\underbrace{\cd{X_k}{X_j}}_{``Z\,"}},X_l) + 2\lambda_k\delta_{kj}\,\gL(\underbrace{\cds{X_i}{T}}_{\lambda_i X_i},X_l)\nonumber\\
&=&\!\! \frac{1}{2}\Big[X_i(\gL(Z,X_l)) + Z(\gL(X_l,X_i)) - X_l(\gL(X_i,Z))\nonumber\\
&&\hspace{.2in}-\gL(Z,[X_i,X_l]) - \gL(X_l,[Z,X_i]) + \gL(X_i,[X_l,Z])\Big]\nonumber\\
&&\hspace{2in}+\, 2\lambda_k\lambda_i\delta_{kj}\delta_{il}\nonumber\\
&=&\!\!g(\cd{X_i}{\!\cd{X_k}{X_j}},X_l) + 2\lambda_i\lambda_k\delta_{kj}\delta_{il}.\label{eqn:term2}
\eeqa
For $(c)$, we observe that because $[X_k,X_i] \in T^{\perp}$, it is of the form $[X_k,X_i] = \sum_{r=1}^{n-1} a_rX_r$ with the $a_r$'s smooth locally defined functions, which in turn implies that $\cds{[X_k,X_i]}{X_j} = \cd{[X_k,X_i]}{X_j} +2a_j\lambda_j T$.  Hence
\beqa
\label{eqn:term3}
\gL(\cds{[X_k,X_i]}{X_j},X_l) \!\!&=&\!\! g(\cd{[X_k,X_i]}{X_j},X_l).
\eeqa
Taken together, \eqref{eqn:term1}, \eqref{eqn:term2}, and \eqref{eqn:term3} yield
\beqa
\label{eqn:Ric4}
\text{Rm}_{\scalebox{0.4}{\emph{L}}}(X_k,X_i,X_j,X_l) \!\!&=&\!\! \text{Rm}(X_k,X_i,X_j,X_l) + 2\lambda_i\lambda_k(\delta_{ij}\delta_{kl}-\delta_{kj}\delta_{il})\nonumber\\
&=&\!\! \text{Rm}(X_k,X_i,X_j,X_l) + (\nabla T^{\flat} \kn \nabla T^{\flat})(X_k,X_i,X_j,X_l).\nonumber\\
\eeqa 
\end{enumerate}
Given \eqref{eqn:Ric1}, \eqref{eqn:Ric2}, and \eqref{eqn:Ric4}, and the fact that $(\nabla T^{\flat} \kn \nabla T^{\flat})(T,\cdot,\cdot,\cdot) = 0$, the proof is complete.
\end{proof}

\section{Examples}
\label{sec:ex}
We provide three examples illustrating our Proposition, all in the noncompact setting, beginning with the canonical example:

\vskip 6pt

{\bf Example~1.}~Consider $n$-dimensional \emph{de Sitter spacetime} $(\mathbb{S}_1^n,\gL)$ ($n \geq 3$), which is the warped product $(\RR \times \mathbb{S}^{n-1},-dt^2 + f(t) \mathring{g})$, where  $\mathring{g}$ is the standard (round) Riemannian metric on $\mathbb{S}^{n-1}$ and 
$$
f(t) = r^2 \cosh^2 \left(\frac{t}{r}\right),
$$
with $r$ the radius of $\mathbb{S}^{n-1}$; consult, e.g., \cite{hartman} and \cite[p.~183ff.]{beem}. This is a geodesically complete Lorentzian manifold with constant positive curvature ($\lambda = 1/r^2$); furthermore,
$$
T \defeq \cds{}{t} = -\partial_t
$$
is a unit timelike gradient, hence has geodesic flow and integrable normal bundle.  The corresponding Riemannian metric 
$$
g \defeq \gL + 2T^{\flat_{\scriptscriptstyle L}} \otimes T^{\flat_{\scriptscriptstyle L}} = dt^2 + f(t)\mathring{g}
$$
is also geodesically complete (see \cite[Lemma~40, p.~209]{o1983}).  In fact the eigenvalues $\lambda_i$ of $\mathcal{D}_{\scalebox{0.4}{\emph{L}}}$ are easily determined here; we shall compute them in Section \ref{sec:conc}.  Finally, note that on any warped product $(\RR \times N,-dt^2 + f(t)h)$, the gradient $T \defeq \cds{}{t} = -\partial_t$ is unit timelike. 

\vskip 6pt

{\bf Example~2.}~On $\RR^3 = \{(x^1,x^2,x^3)\}$, consider a Lorentzian metric $\gL$ whose components $({\gL}_{ij})$ with respect to the coordinate basis $\{\partial_1,\partial_2,\partial_3\}$ take the form
$$
({\gL}_{ij}) = \begin{pmatrix}
0 & 1 & 0\\
1 & H(x^1) & 0\\
0 & 0 & H(x^1)/2
\end{pmatrix},
$$
where $H(x^1)$ is a smooth function, to be determined below, such that $\gL$ will have constant positive curvature on an open subset of $\RR^3$.  In dimension 3, this is equivalent to being an Einstein metric with positive Einstein constant, which we can take to be $1$:
\beqa
\label{eqn:Ein3}
\text{Ric}_{\scriptscriptstyle L} = \gL.
\eeqa
Though we forego the computations here, it is straightforward to show that \eqref{eqn:Ein3} will hold if and only if
$$
H(x^1) = (x^1+a)^2/2 \comma a \in \RR.
$$
For example, if we take $a = 2$, then $\gL$ will be a Lorentzian metric in the open subset $\{(x^1,x^2,x^3) : x^1 > -2\}$.  Finally, let $f(x^1)$ be a smooth function and consider its $\gL$-gradient:
$$
\nabla^{\scalebox{0.4}{\emph{L}}}\!f = -\frac{f'}{2}(x^1+2)^2\,\partial_1 + f'\,\partial_2.
$$
Another computation shows that $\gL(\nabla^{\scalebox{0.4}{\emph{L}}}\!f,\nabla^{\scalebox{0.4}{\emph{L}}}\!f) = -1$ if and only if
$$
(f')^2 = \frac{2}{(x^1+2)^2}\cdot
$$
Taking the smooth solution $f(x^1) = \sqrt{2}\,\text{ln}\big(\frac{x^1+2}{2}\big)$, we thus have that the pair $(\gL, \nabla^{\scalebox{0.4}{\emph{L}}}\!f)$, when restricted to $\{(x^1,x^2,x^3) : x^1 > -2\}$, yields a Lorentzian manifold with constant positive curvature and a unit timelike vector field with geodesic flow and integrable normal bundle.  The corresponding Riemannian metric $g$ is then
$$
g \defeq \gL + 2df \otimes df.
$$
With respect to $g$, $\nabla^{\scalebox{0.4}{\emph{L}}}\!f$ will also be a unit length vector field with geodesic flow and integrable normal bundle.

\vskip 6pt

{\bf Example~3.}~Our last example is not one of constant curvature, but it does involve an important class of Lorentzian manifolds.  On $\RR^4 = \{(v,u,x,y)\}$, consider the Lorentzian \emph{pp-wave} metrics
\beqa
\label{eqn:pp}
\gL = H(u,x,y) du \otimes du + du \otimes dv + dv \otimes du + dx \otimes dx + dy \otimes dy,
\eeqa
where $H(u,x,y)$ is a smooth function independent of $v$.  For the role such metrics play in modeling gravitational waves, as well as their geometric properties more generally, consult \cite{AMS} and \cite{flores}.  Define the vector field
$$
T \defeq \frac{1}{2}(H +f^2 +h^2 +1)\partial_v - \partial_u +f\partial_x + h\partial_y,
$$ 
with $f(u,x,y)$ and $h(u,x,y)$ smooth functions, to be determined below, such that $T$ will have geodesic flow and integrable normal bundle (note that  $\gL(T,T) = -1$ for any $f,h$). We now sketch the conditions on $f, h,$ and $H$ needed for these two properties to hold.  To begin with, the normal bundle $T^{\perp_{\scriptscriptstyle L}}$ is spanned by the orthonormal frame $\{X,Y,Z\}$ given by
\beqa
X = f\partial_v + \partial_x &\comma& Y = h\partial_v+\partial_y,\nonumber\\
&&\hspace{-1.85in} Z = \frac{1}{2}(H +f^2 +h^2 -1)\partial_v - \partial_u +f\partial_x + h\partial_y.\nonumber
\eeqa
(In general, these vector fields are not eigenvectors of $\mathcal{D}_{\scalebox{0.4}{\emph{L}}}$ in \eqref{eqn:D0}; e.g., $\cds{X}{T} = f_xX+h_xY$ and $\cds{Y}{T} = f_yX + h_yY$.) To ensure integrability of $T^{\perp_{\scriptscriptstyle L}}$ via Frobenius's Theorem, we seek to satisfy the Lie bracket conditions
$$
\gL(T,[X,Y]) = \gL(T,[X,Z]) = \gL(T,[Y,Z]) = 0. 
$$
Computing these yields (consult, e.g., \cite[Chapter 13]{beem} for all relevant covariant derivatives), respectively, the following necessary and sufficient conditions on $f, h,$ and $H$:
\beqa
\label{eqn:big3}
h_x = f_y \comma ff_x + hf_y - f_u = \frac{H_x}{2} \comma hh_y + fh_x - h_u = \frac{H_y}{2}\cdot
\eeqa
Incidentally, the latter two equations also ensure that $T$ has geodesic flow: $\cds{T}{T} = 0$.  Many functions $f,h$, and $H$ exist satisfying \eqref{eqn:big3}; e.g., if $f,h,$ and $H$ are functions of $x$ and $y$ only, then we may take any $f,h$ satisfying $h_x = f_y$ and set $H = f^2 +h^2$.  One particularly interesting case\,---\,a so called \emph{plane wave}\,---\,is given by the choices
$$
f(x,y) = y \comma h(x,y) = x \comma H(x,y) = x^2 + y^2,
$$
because, not only will $T$ have geodesic flow and integrable normal bundle, but for this choice of $H$ the metric $\gL$ is in fact geodesically complete and has vanishing \emph{$N$-Bakry-\'Emery tensor} $\text{Ric}^{7/2}_u$ (i.e., with  ``synthetic dimension" $N = 7/2$):
$$
\text{Ric}^{7/2}_u \defeq \text{Ric}_{\scalebox{0.4}{\emph{L}}} + \text{Hess}\,u - \frac{du \otimes du}{7/2 - 4} = 0.
$$
(For the definition and properties of such tensors, see \cite{WW}.  In fact $\nabla^{\scalebox{0.4}{\emph{L}}} u = \partial_v$ is a parallel ``lightlike" vector field for all pp-waves \eqref{eqn:pp}, so that its Hessian, $\text{Hess}\,u$, vanishes.) The corresponding Riemannian metric
$$
g \defeq \gL + 2T^{\flat_{\scriptscriptstyle L}} \otimes T^{\flat_{\scriptscriptstyle L}}
$$
satisfies $\text{Ric} = \text{Ric}_{\scriptscriptstyle L}$, though the Riemannian gradient $\nabla u$ is not parallel.

\section{Proof of Theorem}
\label{sec:proof}
\begin{proof}[\textcolor{blue}{Proof of Theorem.}]
Suppose a pair $(g,T)$ satisfying \eqref{eqn:1} exists on a compact manifold $M$ of dimension $\geq 3$.  We first consider the case $\lambda \leq 0$; it turns out that the obstruction in this case occurs at the level of the Ricci tensor, via a well known equation (\eqref{eqn:Bochner1} below) and Riccati analysis, as follows.  Set $\gL \defeq g - 2T^{\flat}\otimes T^{\flat}$.  Then
\beqa
\text{Ric}(T,T) &=& \sum_{i=1}^{n-1} \text{Rm}(X_i,T,T,X_i)\nonumber\\
&\overset{\eqref{eqn:Ric1}}{=}& \sum_{i=1}^{n-1} \text{Rm}_{\scalebox{0.4}{\emph{L}}}(X_i,T,T,X_i) = \text{Ric}_{\scriptscriptstyle L}(T,T),
\eeqa
so that if $(M,\gL)$ has constant curvature $\lambda \leq 0$\,---\,more generally, if $(M,\gL)$ is Einstein with nonpositive Einstein constant\,---\,then
\beqa
\label{eqn:Ric00}
\text{Ric}(T,T) = \text{Ric}_{\scriptscriptstyle L}(T,T) = (n-1)\lambda \underbrace{\,\gL(T,T)\,}_{-1} \geq 0.
\eeqa
Next, setting $i=j$ in \eqref{eqn:Bochner} and summing over $i = 1,\dots,n-1$ yields the following Bochner-type equation,
\beqa
\label{eqn:Bochner1}
T(\text{div}\,T) = -\text{Ric}(T,T) - \sum_{i=1}^{n-1} \lambda_i^2,
\eeqa
where we've used the fact that 
$$
\text{div}\,T = \sum_{i=1}^{n-1} \lambda_i = \text{tr}_g \nabla T^{\flat}.
$$
Now, via the Schwarz inequality $$\sum_{i=1}^{n-1} \lambda_i^2 \geq \frac{1}{n-1}(\lambda_1+ \cdots +\lambda_{n-1})^2,$$ \eqref{eqn:Bochner1} reduces to
\beqa
\label{eqn:Bochner3}
\kk(\text{div}\,\kk)\ \leq \ -\text{Ric}(\kk,\kk) - \frac{(\text{div}\,\kk)^2}{n-1},
\eeqa
which permits, in turn, the following well known Riccati analysis: since $T$ is complete ($M$ being compact),
$$
\text{Ric}(T,T) \geq 0 \overset{\eqref{eqn:Bochner3}}{\imp} \text{div}\,T = 0 \overset{\eqref{eqn:Bochner1}}{\imp} \text{Ric}(T,T) = \lambda_i = 0.
$$
It follows that if $\text{Ric}(T,T) = 0$, then $T$ must be parallel, and that the case $\text{Ric}(T,T) > 0$ cannot occur.  The latter implies that $\lambda < 0$ cannot occur; the  former, that if $\lambda = 0$, then $\nabla T^{\flat}$ is zero and \eqref{eqn:1} vanishes, so that $(M,g)$ is flat.  Furthermore, its universal covering splits isometrically as a product $\RR \times N$, by the de Rham Decomposition Theorem (see \cite[p.~384]{Petersen}).  This settles the case $\lambda \leq 0$.  For the case $\lambda > 0$, we will employ a different strategy; indeed, since $\text{Ric}(T,T) < 0$ when $\lambda > 0$ (via \eqref{eqn:Ric00}), the Riccati analysis we applied to \eqref{eqn:Bochner3} is unavailable here.  Instead, we substitute $g = \gL + 2T^{\flat}\otimes T^{\flat}$ into our Proposition (note that $T^{\flat_{\scriptscriptstyle L}}\otimes T^{\flat_{\scriptscriptstyle L}} = T^{\flat}\otimes T^{\flat}$), to obtain
\beqa
\text{Rm}_{\scalebox{0.4}{\emph{L}}} \!\!&=&\!\! \text{Rm} + \nabla T^{\flat} \kn \nabla T^{\flat}\nonumber\\
&\overset{ \eqref{eqn:1}}{=}&\!\!\frac{1}{2}\lambda g \kn g - 2\lambda g \kn (T^{\flat}\otimes T^{\flat})\nonumber\\
&=&\!\!\frac{1}{2}\lambda \gL \kn \gL + 2\lambda \gL \kn (T^{\flat}\otimes T^{\flat}) - 2\lambda \gL \kn (T^{\flat}\otimes T^{\flat})\nonumber\\
&=&\!\!\frac{1}{2}\lambda \gL \kn \gL,\nonumber
\eeqa
where we've used the fact that $(T^{\flat}\otimes T^{\flat}) \kn (T^{\flat}\otimes T^{\flat}) = 0$.  But as mentioned in the Introduction, such a (compact) Lorentzian manifold is impossible when $\lambda > 0$, by \cite{CM} and \cite{Klingler}.
\end{proof}

\section{Concluding Remarks}
\label{sec:conc}

\vskip 6pt
{\bf Remark~1.}~Observe that for a curvature 4-tensor of the form \eqref{eqn:1}, the term $g \kn (T^{\flat}\otimes T^{\flat})$ is zero on the frame $\{T,X_1,\dots,X_{n-1}\}$ except for the component
$$
g \kn (T^{\flat}\otimes T^{\flat})(X_i,T,T,X_i) = 1.
$$
In particular, the sectional curvature of any 2-plane containing $T$ is 
$$
\lambda - 2\lambda = -\lambda,
$$
as opposed to $\lambda$, which would have been the case with constant curvature $\text{Rm} = \frac{1}{2}\lambda g \kn g$.  On the other hand, the term $\nabla T^{\flat} \kn \nabla T^{\flat}$ is zero except for the component
$$
(\nabla T^{\flat} \kn \nabla T^{\flat})(X_i,X_j,X_j,X_i) = 2\lambda_i\lambda_j \comma i \neq j,
$$
so that 2-planes spanned by $\{X_i,X_j\}$ now have sectional curvature
$$
\lambda - 2\lambda_i\lambda_j.
$$
For this reason, we may regard \eqref{eqn:1} as a curvature tensor for which constant curvature has been ``broken by $T$."

\vskip 6pt

{\bf Remark~2.}~Our Theorem is uninteresting in dimension 2.  Indeed, although any Riemannian 2-manifold satisfies $\text{Rm} = \frac{1}{2}K g \kn g$ with $K$ the Gaussian curvature (see, e.g., \cite[p.~250]{Lee}), and although $K$ is neither constant nor signed in general, nevertheless $K = -\lambda$ if \eqref{eqn:1} is satisfied, because in dimension 2, $$
\nabla T^{\flat} \kn \nabla T^{\flat} = 0 \comma  g \kn (T^{\flat} \otimes T^{\flat}) = \frac{1}{2}g \kn g.
$$
But by the Gauss-Bonnet Theorem, when the Euler characteristic is zero (i.e., the manifold is a 2-torus or a Klein bottle), the only compact Riemannian 2-manifold with constant curvature is the flat one.

\vskip 6pt

{\bf Remark~3.}~Via \eqref{eqn:rel1}, the Ricci tensor of $g$ is related to that of $\gL$ as follows:
\beqa
\label{eqn:RicRic}
\text{Ric} = \text{Ric}_{\scalebox{0.4}{\emph{L}}} + 2\text{Rm}_{\scalebox{0.4}{\emph{L}}}(T,\cdot,\cdot,T) - \text{tr}_g(\nabla T^{\flat} \kn \nabla T^{\flat}).
\eeqa
It is not signed in general.  Indeed, in the specific case of \eqref{eqn:1} with $\lambda > 0$, note that while $\text{Ric}(T,T) < 0$ by \eqref{eqn:Ric00} (because now $\lambda > 0$), $\text{Ric}(X_i,X_i)$ is given by
\beqa
\text{Ric}(X_i,X_i) \!\!&=&\!\! \underbrace{\,\text{Rm}(T,X_i,X_i,T)\,}_{\text{$\text{Rm}_{\scalebox{0.4}{\emph{L}}}(T,X_i,X_i,T)$}} + \sum_{k\neq i} \underbrace{\,\text{Rm}(X_k,X_i,X_i,X_k)\,}_{\text{$\text{Rm}_{\scalebox{0.4}{\emph{L}}}(X_k,X_i,X_i,X_k)\,-\,2\lambda_i\lambda_k$}}\nonumber\\
&=&\!\! -\lambda + (n-2)\lambda - 2\lambda_i\sum_{k\neq i} \lambda_k,\label{eqn:sign}
\eeqa
where we observe that, because $\gL(T,T) = -1$,
$$
\text{Rm}_{\scalebox{0.4}{\emph{L}}}(T,X_i,X_i,T) = \frac{1}{2}\lambda (g \kn g)(T,X_i,X_i,T) = -\lambda.
$$
In particular, \eqref{eqn:sign} need not be negative in general, so that the Ricci tensor need not be signed.  In fact this is evident in Example 1 of Section \ref{sec:ex} above, namely, de Sitter spacetime $(\RR \times \mathbb{S}^{n-1},-dt^2+f(t)\mathring{g})$, here with radius $r = 1$ so that $\lambda = 1$.  By standard formulae for covariant derivatives on warped products (see, e.g., \cite[p.~206]{o1983}), for any $X \in T^{\perp_{\scriptscriptstyle L}}$,
$$
\cds{X}{T} = \frac{T(f)}{f}X = -2\tanh(t) X,
$$
so that each $\lambda_i = -2\tanh(t)$, and hence
$$
\text{Ric}(X_i,X_i) = (n-3) - 8(n-2)\tanh^2(t).
$$
For $n \geq 4$, this can be positive, e.g., at $t = 0$.  In general, when \eqref{eqn:1} holds, \eqref{eqn:RicRic} takes the form
$$
\text{Ric} = (n-1)\lambda g + \lambda (g \kn g)(T,\cdot,\cdot,T) - \text{tr}_g(\nabla T^{\flat} \kn \nabla T^{\flat}).
$$

\vskip 6pt

{\bf Remark~4.}~Regarding \eqref{eqn:Bochner3}, note that if $\lambda > 0$, then even in the case of equality,
$$
\kk(\text{div}\,\kk) =  (n-1)\lambda - \frac{(\text{div}\,\kk)^2}{n-1},
$$
nontrivial complete solutions exist.  Indeed, if $s$ is an affine parameter along an integral curve of $T$, then 
$$
(n-1)\sqrt{\lambda}\,\text{tanh}\big(\sqrt{\lambda}\,s+c\big) \comma \pm(n-1)\sqrt{\lambda}
$$
are complete solutions, so that the Riccati analysis above is unavailable here.
\bibliographystyle{alpha}
\bibliography{Riemannian_n}
\end{document}